\documentstyle[12pt]{article}

\title{Braided modules and reflection equations}
\author{D.~Gurevich\\
ISTV, Universit\'e de Valenciennes\\
59304 Valenciennes, France}
\date{}

\begin{document}

\newtheorem{proposition}{Proposition}
\newtheorem{definition}{Definition}
\newtheorem{conjecture}{Conjecture}
\newtheorem{remark}{Remark}
\newcommand{\ren}{\rho^{End}_q}
\newcommand{\orr}{\overline{\rho}}
\newcommand{\br}{[\,\,,\,\,]}
\newcommand{\brq}{[\,\,,\,\,]_q}
\newcommand{\vqp}{V^q_+}
\newcommand{\vqm}{V^q_-}
\newcommand{\gggg}{ g}
\newcommand{\tS}{\widetilde{S}}
\newcommand{\ahq}{A_{h,q}}
\newcommand{\ahc}{A_{h,1}^c}
\newcommand{\rn}{\rho_{\nu}}
\newcommand{\usl}{U_q(sl(2))}
\newcommand{\ug}{U(\gggg)}
\newcommand{\us}{U(sl(2))}
\newcommand{\vv}{V^{\otimes 2}}
\newcommand{\ww}{W^{\otimes 2}}
\newcommand{\uqs}{U_q(sl(2))}
\newcommand{\ork}{\overline{\rho_k}}
\newcommand{\uog}{U(\overline{\gggg})}
\newcommand{\uos}{U(\overline{sl(2)})}
\newcommand{\uq}{U_q(\gggg)}
\newcommand{\sn}{sl(n)}
\newcommand{\ogg}{\overline{\gggg}}
\newcommand{\osl}{\overline{sl(2)}}
\newcommand{\osln}{\overline{sl(n)}}
\newcommand{\rd}{\rho^{\otimes 2}}
\newcommand{\vbq}{V_{\beta}^q}
\newcommand{\uqq}{U^q_{1/2}}
\newcommand{\uqk}{U^q_{k/2}}
\newcommand{\vq}{V^q}
\newcommand{\vb}{V_{\beta}}
\newcommand{\vqb}{V_{\beta}^q}
\newcommand{\voq}{V_{\omega}^q}
\newcommand{\vkoq}{V_{k\omega}^q}
\newcommand{\vkb}{V_{k\beta}}
\newcommand{\vbk}{V_{\beta}^{\otimes k}}
\newcommand{\vko}{V_{k\omega}}
\newcommand{\vlo}{V_{l\omega}}
\newcommand{\vklo}{V_{(k+l)\omega}}
\newcommand{\aq}{A_{0,q}}
\def\ot{\otimes}
\def\De{\Delta}
\def\qq{q^{-1}}

\maketitle

\begin{abstract}
We introduce a representation theory 
of q-Lie algebras defined earlier in \cite{DG1}, \cite{DG2}, 
formulated in terms of braided modules.
We also discuss  other ways to define Lie algebra-like objects
 related to quantum groups, in particular, those based on 
the so-called reflection
equations. We also investigate the truncated tensor product of 
braided modules. 

\end{abstract}

\section{Introduction}
Since the creation of the quantum groups theory enormous 
efforts have been made  to find
 proper generalizations of Lie algebra structures  
(in spirit of super-theory)
connected to the deformed enveloping algebras $\uq$. More precisely, the
following problem was considered: find an 
operator $\br :\vv\to V$ satisfying
a twisted (quantum) version of Jacobi identity,
which is somehow related  with
the algebras $\uq$. Recently some hopes to find a consistent construction
of such Lie algebra-like objects were connected with the so-called reflection 
equations (RE).

In particular, these equations participate  in the construction of 
the bicovariant
differential calculus initiated  by S.~Woronowicz (cf. f.e. \cite{IP1},
\cite{IP2} and the references therein).  The braided groups of S.~Majid can 
 also be
treated in terms of the RE. More precisely, the braided groups can be 
represented as certain quotient algebras of graded quadratic algebras 
defined by the
RE (RE algebras for short). There exists a natural but a somewhat
 trivial way to
assign   certain
quadratic-linear algebras to these graded quadratic ones,
 which can be thought of as enveloping algebras of some Lie
algebra-like objects (we represent here a version of Jacobi identity 
which is
valid for  objects of such type).  However, though the RE algebras represent 
interesting (and, hopefully, flat in the $sl(n)$-case) deformations  of the
function algebra on the group $SL(n),$ their further transformation into
the mentioned quadratic-linear algebras is not meaningful from
the deformation theory point of view (cf. Section 5).

Another way to construct q-deformed Lie algebra-like objects called
{\em braided (or q-)Lie algebras} was suggested in \cite{DG2} 
(the $sl(2)$-case
was previously considered in \cite{DG1}). We present here a
construction of braided Lie algebras and make an attempt to develop
their representation theory. The linear spaces, where our braided Lie
algebras are represented, are called {\em braided modules}. Let us note that
the family of  braided modules does not form any tensor category and
is a subset in the category of  $\uq$-modules (apart from the
$sl(2)$-case, when any $U_q(sl(2))$-module can be equipped 
with a {\em braided
structure}).

We want to emphasize that our approach 
to the definition of q-analogues of
ordinary Lie  algebras is in fact the most natural way to do it
(at least for objects "living" in the category of $\uq$-modules)
although any form of Jacobi identity looking like the classical one does not
exist for them (apart from the $sl(2)$-case). The reason for this is the
following: the $\gggg$-module $End(\gggg)$ contains "too many" irreducible
$\gggg$-components (cf. Remark 1 where a truncated version of Jacobi identity is
discussed).

The paper is organized as follows. In Section 2 we present a construction of
braided Lie algebras and compare it with some other approaches to the Lie
algebra-like objects. In Section 3 we discuss a representation theory 
for braided
Lie algebras in terms of braided modules. In Section 4 we consider some
truncated version of the tensor product for them. We complete 
the  paper with
a discussion on the  RE algebras.

Throughout the paper the ground field $k$ is assumed to be $R$ or $C$ and the
parameter $q$ is generic. The category of $\uq$-modules is
 denoted by $\uq$-Mod.

\section{Braided Lie algebras}

Let us recall the construction of braided (or q-) Lie brackets introduced in
\cite{DG1} and \cite{DG2}.

Let $\gggg$ be a simple Lie algebra. Let us assume that
it is equipped with the left adjoint representation 
$\rho(a)z=[a,z],\, a,\,z\in
\gggg$ as an object of the category of
$\gggg$-modules. This object will be denoted by $V$. Let us consider the
extension $\rho : U(\gggg) \to End(V)$ of this representation to
the enveloping algebra and introduce a representation
$\rho_q: \uq \to End(V)$ deforming the last one.

For the classical (non-exceptional) Lie algebras this deformation can be
constructed as follows. Fix a vector fundamental representation
$\gggg\to End(U)$ of $\gggg$. Then the space $End(U)$ can be also equipped
with a  $\gggg$-module structure. Then one of the components  
the decomposition of $End(U)$ into a direct sum of irreducible 
$\gggg$-modules $End(V)=
\oplus V_{\beta}$ is isomorphic to $\gggg$ itself. Hereafter we denote 
 the $\gggg$-module of highest weight $\beta$ by
$V_{\beta}$.

Now let us equip the space $U$ with a $\uq$-module structure 
(it is well known
that the action of the generators $H_{\alpha},\,X_{\pm\alpha}\in\uq$ on the
fundamental module coincides with the classical one). 
Let    $\rho_q:
\uq\to End(U)$ denote this representation.  Then, using the coproduct
operator and the antipod of the algebra $\uq$, we equip the space
$End(U)$ with a $\uq$-module structure as follows:

Introduce a representation $\ren: \uq\to End(End(U))$ putting
$$\ren(a)M=\rho(a_1)\circ M\circ\rho(\gamma(a_2)),\,a\in\uq, \,M\in End(U).$$
We denote  the matrix product by  $\circ$, $\gamma$ is the
antipode in $\uq$ and $a_1\ot a_2$ is the Sweedler's notation for $\De(a)$.

Moreover, we deal with a
coordinate representation of module elements. We consider the endomorphisms
as matrices and their action is the left multiplication by these matrices.

Let us remark that this way to equip $End(U)$ with a structure 
of $\uq$-module
is compatible with the matrix product in it in the following sense:
$$\ren(a)(M_1\circ M_2)=\ren(a_1)M_1\circ\ren(a_2)M_2.$$
This means that the product $\circ: End(U)^{\ot 2}\to End(U)$ is
a $\uq$-morphism.

Thus, we have turned the space $End(U)$ into a $\uq$-module. 
(Similarly
we can equip  the space $End(U)$ with a $\uq$-module structure for any
$\uq$-module $U$.)  Decompose the space $End(U)=\oplus\vbq$ into a
direct sum of $\uq$-modules. One of this component can be identified with $V$
and it is equipped with a $\uq$-module structure.

For exceptional Lie algebras this method to deform the adjoint representation
does not work, but anyway for any simple Lie algebra $\gggg$ 
the adjoint module
$V$ can be equipped with the structure of a $\uq$-module. Using the coproduct
operator in $\uq$, extend the $\uq$-module structure to the space $\vv$.
Now decompose this $\uq$-module into a direct sum of irreducible modules
$V_{\beta}^q$. Then the space $V$ itself enters this decomposition with
multiplicity 2 (if $\gggg=sl(n),\, n>2$) or 1 (otherwise). 
However, in the first
case one of the modules, isomorphic to $V$, belongs to $I_+^q$ 
and the other 
belongs to $I_-^q$, where $I_+^q \,(I_-^q)$ denotes the q-analogue of 
symmetric
(skew-symmetric) subspace in $\vv$ (for definitions of these
 subspaces see, e.g., \cite{DG2}). Let us denote the mentioned 
$\uq$-modules by
$\vqp$ and $\vqm$, respectively (the notation $\vqm$ will be kept also in the
case when the module is unique).

Then it is natural to introduce the {\em braided Lie bracket} as a
$\uq$-morphism $\brq:\vv\to V$ killing all $V^q_{\beta}$ 
excluding that $\vqm$
(this property defines the bracket up to a non-trivial factor).

It is evident that this bracket is $\tS$-skew-symmetric 
(i.e. $\brq=-\brq\tS$)
where $\tS$ is the operator defined as follows:
$\tS=id$ on $I_+^q$ and $\tS=-id$ on $I_-^q$.

The space $V$ equipped with the bracket $\brq$ will be denoted by $\ogg$.

Now let us introduce  the enveloping algebra $\uog$ of
the braided Lie algebra $\ogg$. Put
$$\uog=T(\ogg)/\{\vbq-\tau\brq \vbq,\,\,\forall\,\,\vbq\subset I^q_-\}.$$
Hereafter we let $T(V)$ denote the free tensor algebra of the space $V$ 
and 
$\{J\}$ denote the ideal  generated by the set $J\subset T(V)$.
If to compare this definition to a similar definition in 
\cite{DG1} and \cite{DG2}, one can see that we have
introduced a factor $\tau$ in this definition. 
It will be discussed in Section 4
(cf. Remark 1).

Let us point out two properties of the algebra $\uog$. 
First, it is evident that
the product in it is $\uq$-invariant, i.e.,
$a(x\circ y)=a_1(x)\circ a_2(y)$ for any $a\in\uq,\,\,x,y\in \uog$. 
Second, these
algebras are not the result of any flat deformation (the reader is referred 
to
\cite{DG2} for definitions) of the ordinary enveloping algebra since even
the deformation of the symmetric algebra of the space $V$ to its q-analogue
$T(V)/\{I_-^q\}$ is not flat (apart from the $sl(2)$-case),
 but now it does not
matter for us. (However, certain quotient algebras of the algebra
$\uog ,$ corresponding to the $R-$matrix type orbits 
in terminology of \cite{GP} , can be
in principle flatly deformed to their q-counterparts. 
The problem is to find a
proper description of such "quantum orbits" in spirit of the paper \cite{DG2}.
This problem will be considered in details in other publications.)

Let us remark that there exists a number of papers (cf. \cite{M}, \cite{DH},
\cite{LS} and others)  where Lie algebra-like objects arise
as some subspaces of the quantum group $\uq$. More precisely, the approach
suggested in these papers involves searching for a finite dimensional space
$L\subset\uq$ invariant with respect to the q-adjoint operator in $\uq$
(sometimes certain complementary conditions on $L$ are imposed, see, e.g.,
\cite{DH}). Then this space equipped with the bracket defined by the 
q-adjoint
operator is regarded as a "quantum Lie algebra". 

Apparently, the objects arising from this approach are close to ours. 
 However,
in our approach the space $V$ is not embedded in the 
quantum group.
We use $\uq$ only for a suitable description of the objects and 
morphisms of the
category $\uq$-Mod\footnote{Let us remark that our method 
can be applied to define a q-Lie bracket in some non-quasiclassical category,
i.e., those generated by a space  $V$ equipped with a non-deformational 
quantum
R-matrix, say, that of Hecke type, possessing a non-standard  
Poincar\'e series
(such quantum matrices have been constructed in \cite{G}). It would be
interesting to find a q-deformation of generalized Lie brackets computed
explicitly in \cite{G} for some involutive solutions of the quantum 
Yang-Baxter
equation.}.

Let us also remark that in contrast with the approach developped, say, in
\cite{M}, where certain Lie algebra-like objects (called quantum and braided Lie
algebras)  are introduced by means of a version of Jacobi identity our braided
Lie algebras are introduced directly, by constraction. They satisfy in general 
a very truncated form (compared with the classical one) of Jacobi identity
(cf. Remark 1).

\section{Braided modules}

In the present Section we discuss a representation theory of the braided
algebra $\ogg$ (essentially for $\gggg=sl(n)$). Compared with quantum groups,
 the braided algebras look more
like super- or generalized (S-)Lie algebras defined in 
\cite{G} for an involutive solution $S$ to the QYBE. However, unlike
the latter objects the braided Lie algebras admit only a truncated version of
the representation theory. Our goal is to equip $\uq$-modules with the
 structure
of $\uog$-modules but only certain of them can be "converted" 
into the braided ones.

\begin{definition} We say that an $\uq$-module $U^q$ is a braided module or,
more precisely, a braided $\ogg$-module if it can be equipped with 
the structure of a
$\uog$-module in such a way that the corresponding map
$\uog\to End(U^q)$ is $\uq$-morphism. 
We also say  that the classical counterpart $U=U^1$ of the
$\uq$-module $U^q$ allows braiding.
\end{definition}

A natural way to construct braided modules is given by the following

\begin{proposition} Let $U^q$ be a $\uq$-module. If the decomposition
$End(U^q)=\oplus V^q_{\gamma}$ of the $\uq$-module $End(U^q)$ into a 
direct sum
of irreducible $\uq$-modules is such that\\
1. it does not contain the modules isomorphic to $\vbq\subset I_-^q$
apart from
those isomorphic to $V^q$,\\
2. the multiplicity of the module $V^q$ is 1,\\
then $U^q$ can be equipped with a braided module structure (briefly, 
braided structure).
\end{proposition}

{\bf Proof.} By the assumption  there exists a unique
$\uq$-module in  $End(U^q),$ isomorphic to $V^q$. Consider an 
$\uq$-morphism defined up
to a factor
\begin{equation}
\rho: V^q\to End(U^q).
\end{equation}
This map is an almost representation of the braided Lie
algebra $\ogg$ in the sense of the following

\begin{definition} We say that a map (1) is an almost representation of the
braided Lie algebra $\ogg$ if it is a $\uq$-morphism and the following
properties are satisfied\\
1. $\circ\rd(\vbq)=0$ for all $\vbq\subset I_-^q$ apart from that
$\vbq=V^q_-$,\\
2. $\circ\rd V_-^q=\nu\rho[\,\,,\,\,]_qV^q_-$ with some $\nu\not=0$\\
(if $\nu=\tau,$ where $\tau$ will be defined below, then 
we have a representation).
\end{definition}

In a more explicit form these conditions can be reformulated as follows. 
If the
elements $\{b^{i,j}_{k,\beta}u_iu_j, 1\leq k\leq dim\, \vbq\}$ form a 
basis of
the space $\vbq\subset I_-^q$ and analogously the
elements $\{b^{i,j}_{k,-}u_iu_j, 1\leq k\leq dim\, \vbq\}$ form a 
basis of the space $V_-^q$ then
\begin{equation}
b^{i,j}_{k,\beta}\rho(u_i)\rho(u_j)=0\,\,{\rm if}\, \,\vbq\not=V_-^q\, \,{\rm
and}\,\, b^{i,j}_{k,-}\rho(u_i)\rho(u_j)=\nu b^{i,j}_{k,-}\rho([u_i,u_j]_q).
\end{equation}

Let us complete  the proof. The image of the composed map
$\circ\rd(V^q_-)$ is isomorphic to the $\uq$-module $V^q$ since $\rho$ and
$\circ$ are $\uq$-morphisms. Such a module in the decomposition of
$End(U^q)$ is unique therefore the image of the space $[\,\,,\,\,]_q
V^q_-=V^q$ with respect to the morphism $\rho$ coincides with the previous
 one 
(it suffices to show that the images of the highest weight element of the 
module
$V_-^q$ with respect to both operators coincide up to a factor). 
This gives the
second property of Definition 2 (the property that $\nu\not=0$ for 
a generic $q$
follows from the fact that it is so for $q=1,$ since in this case 
 we
have an ordinary representation $\gggg\to End(U)$ by construction).

The first property follows from the fact that $End(U)$ does not contain 
any modules
isomorphic to $\vbq\subset I^q_-,\, \vbq\not=V_-^q$.  Finally, changing the
scale, i.e., considering the map $\rho_{\nu}=\tau\nu^{-1}\rho$ instead of $\rho$
we get a representation of the braided Lie algebra under consideration. This
completes the proof.

A natural question arises: how many $\uq$-modules can be converted 
into braided
ones or in other words how many $\gggg$-modules allow braiding?
 An answer to this question for $sl(n)$-case is given by the following
Proposition which can be proved by straightforward computations 
using the Young
diagrams technique.

\begin{proposition} Let $\omega$ be a fundamental weight of the Lie
algebra $sl(n)$. Then the $sl(n)$-modules $\vko$ (for any $k\in N$) allow
braiding (hereafter we denote  the family of non-negative integers by $N$).
In other words, their q-analogues $\vko^q$ are braided modules.
\end{proposition}

As for other simple Lie algebras $\gggg$
it seems very plausible that a similar statement is valid 
for  fundamental weights $\omega$ such that their orbits in $\gggg^*$ are of
R-matrix type (in the $sl(n)$-case all fundamental weights satisfy this
condition).

If this is the case, we see that  "good"  modules (i.e,. allowing braiding)
correspond to "good"  orbits (i.e. of R-matrix type) in $\gggg^*$.
A quantum (or braided) analogue of this correspondence may be regarded as 
a quantization of the ordinary orbit method (the $sl(2)$-case has 
been considered
from this point of view in \cite{DGR}).

\section{Truncated product of braided modules: the $sl(2)$-case}

In the present Section we discuss the following problem: 
what is a regular way
to multiply some braided modules? In the classical case, given two
$\gggg$-modules, we can equip their tensor product with 
a $\gggg$-module structure by
means of the ordinary coproduct operator in the enveloping algebra $\ug$.

This procedure cannot be generalized to the braided case because of two
difficulties: first,
the algebra $\uog$ does not have any Hopf structure (ordinary or braided),
and  second, (but in fact it is another interpretation of the first
difficulty) the braided modules do not form any tensor category.
However, there exists a module, 
namely, $\vklo$, in the tensor product $\vko\otimes\vlo$
 (where $\omega$ is a fundamental weight and $k,\,l \in N$),  
which can be equipped with a braided structure (we restrict ourselves to
the $sl(n)$-case).

Of course, the above module can be equipped with such a structure 
by means of
the above procedure but a natural question arises: does there exist a more
direct way to do this? Say, in the classical case, given
 a representation
$\rho :\gggg\to End(V_{\beta})$, we can equip the space $\vkb\subset \vbk$ 
with a
$\gggg$-module structure by putting
\begin{equation}
\rho_k(z)v=kP_k(\rho(z)\otimes id_{k-1})v,\,z\in\gggg,\,v\in\vkb,
\end{equation}
where $P_k: \vbk\to \vkb$ denote the natural projector and let $id_{k}$ 
denote the
identity operator $\vbk\to\vbk$. It is easy to see that the map
$\rho_k :\gggg\to End(\vkb)$ defines a representation of the Lie algebra
$\gggg$.

It is natural to try to generalize this construction to the braided case
assuming $\beta=\omega$ to be a fundamental weight (in this case  both
$\gggg$-modules $V_{\beta}$ and $\vkb$
allow braiding)\footnote{In \cite{DG1}, \cite{DG2} a similar approach was
suggested to define braided vector fields.}.

We call the resulting representation the {\em truncated} (since
it is constructed by means of a projector) tensor product of braided modules.
More generally, we consider here truncated tensor powers
of the module  $\vb$. 

\begin{proposition} Let $\rho: \osln\to End(\voq)$ be an almost 
representation
of the braided Lie algebra $\osln$ and let the projector 
$P_k: (\voq)^{\ot k} \to \vkoq$ be a morphism of the
category $U_q(sl(n))$-Mod. Then the map  $\rho_k$ defined by the formula (3)
is also an almost representation of $\osln$. \end{proposition}

This proposition can be proved in the same way as Proposition 1.

Nevertheless, if the initial map $\rho$ is a representation of the 
braided Lie
algebra $\osln$ we can turn the map $\rho_k$ into a representation 
by a proper
change of the factor $k$. We compute this deformed factor in the 
simplest case
$n=2$. However, first we will consider the above factor $\tau$
in the definition of the enveloping algebra of a braided Lie algebra 
more precisely.

Let us consider the quantum group $\uqs$, i.e., the algebra generated by the 
elements
$H,\,X,\,Y$ satisfying the following relations 
$$[H,X]=2X,\; [H,Y]=-2Y,\; [X,Y]=\frac{q^H-q^{-H}}{q-q^{-1}}.$$
Let us assume that the coproduct in this algebra is of the form
$$\De(X)=X\ot 1+q^{-H}\ot X,\; \De(Y)=1\ot Y +Y\ot q^H,\;
\De(H)=H\ot 1+ 1\ot H.$$

Let us consider a three-dimensional space $V$ with a base $\{u,\,v,\,w\}$
equipped with the following action of the quantum group
$$Hu=2u,\ Hv=0,\ Hw=-2w,\ Xu=0,\ Xv=-(q+\qq) u,\ Xw=v,$$ 
$$Yu=-v,\ Yv=(q+\qq) w,\ Yw=0.$$

Then decomposing the space $\vv$ into a direct sum of irreducible
$\uqs$-modules and introducing the q-Lie bracket as above we can obtain the
following multiplication table for it
$$
[u,u]_q=0,\ [u,v]_q=-q^2Mu,\ [u,w]_q=(q+\qq)^{-1}Mv,$$
$$[v,u]_q=Mu,\ [v,v]_q=(1-q^2)Mv,\ [v,w]_q=-q^2Mw,$$
$$[w,u]_q=-(q+\qq)^{-1}Mv,\ [w,v]_q=Mw,\ [w,w]_q=0,$$
with an arbitrary $M$. 

Then the enveloping algebra of $\osl$ is defined by the following relations
\begin{equation}
q^2u v - v u=- h u,\,(q^3+q)( u w-w u) +(1-q^2)v^2= hv,\, -q^2v w + w v=hw,
\end{equation}
where $h=\tau (q^4+1)M$ (this relation can be deduced from the following
equality: $q^2u v - v u=\tau(q^2[u,v]_q-[v,u]_q)$).

We fix $\tau$ in such a way that the (left) q-adjoint operators
$\rho_q(u)z=[u,z]_q$ define a representation, i.e., the following 
relations hold:
$$ q^2[u,[v,z]_q]_q-[v,[u,z]_q]_q=-\tau[u,z]_q,\, (q^3+q)([u,[w,z]_q]_q-
[w,[u,z]_q]_q)+$$
$$(1-q^2)[v,[v,z]_q]_q=\tau [v,z]_q,\,
-q^2[v,[w,z]_q]_q+[w,[v,z]_q]_q= \tau[w,z]_q,\, \forall z\in\ogg.$$
 This implies the relation $h=M(q^4-q^2+1)$ and therefore we have
$\tau =1-(q^2+q^{-2})^{-1}$.

We treat the latter formulae as a braided version of Jacobi identity. 
This
means that the q-adjoint operators corresponding to the generators 
$u,\,v,\,w$
satisfy the same relations as the generators themselves in the enveloping
algebra  $\uos$.

\begin{remark} As for other simple algebras we can normalize the parameter
$\tau$ requiring that it satisfies the relations
\begin{equation}
b^{i,j}_{k-}[u_i,[u_j,z]_q]_q=\tau b^{i,j}_{k-}[[u_i,u_j]_q,z]_q,
\end{equation}
where $\{b^{i,j}_{k-}u_iu_j, 1\leq k\leq dim\, V^q_-\}$ is a basis of 
the space
$V_-^q$.

Unfortunately, in the general case we are unable to ensure the relations
$$b^{i,j}_{k,\beta}[u_i,[u_j,z]_q]_q=0,\,\,\vbq\not=V_-^q$$ 
(cf. (2)). This is the reason why the "q-adjoint representation" of a braided 
Lie algebra $\ogg\not=\osl$ is not a representation at all. 
Summing up, we can say that the reason of this "pathology" 
is the following: 
the $\gggg$-module $End(\gggg)$ contains "too many" irreducible
$\gggg$-components (apart from the $sl(2)$-case). 

Nevertheless, we treat the formula (5) as a truncated version of Jacobi
identity. The term "truncated" means here that only the space
$V^q_-$ (instead of the whole $I^q_-$) participats in these relations.
 \end{remark}

Let us now discuss  a way to modify the factor $k$ in the formula (3).
Fix a base $\{a,\,b\}$ in a $\usl$-module $\uqq$ such that
$$Xa=0,\,Xb=a,\, Ha=a,\, Hb=-b,\, Ya=b,\, Yb=0.$$
Consider the map $V\otimes \uqq\to \uqq$ defined as follows:
$$u\ot a\to 0,\, u\ot b \to a,\, v\ot a \to q^{-1}a, v\ot b \to -qb,\, 
w\ot a\to
q^{-1}b,\, w\ot b \to 0.$$

It is easy to see that this map is an $\usl$-morphism. It enables us to
construct an almost representation of the braided Lie algebra $\osl$
according to the above procedure. This
almost representation is given by $$
\rho(u)=\left(\begin{array}{cc}0&1\\0&0\end{array}\right),\,
\rho(v)=\left(\begin{array}{cc}q^{-1}&0\\0&-q\end{array}\right),\,
\rho(w)=\left(\begin{array}{cc}0&0\\q^{-1}&0\end{array}\right).$$ 

More precisely, the operators $\rho(u),\, \rho(v),\, \rho(w)$ satisfy the
relations (4) with $h=h_1=q^3+q^{-1}$. Consider now the almost
representation $\ork=P_k(\rho\ot id_{k-1})$, where $P_k$ is the above 
projector
in the category $\uq$-Mod. The operators $\ork(u),\, \ork(v),\, \ork(w)$
satisfy the relations (4) with a factor $h=h_k$. Let us compute it. Then the
proper q-analogue of the factor $k$ in the formula (3) is $h_1/h_k$.

Fix  the base generated by the
elements $f_1=a^k,\,f_2==Ya^k,\,..., f_{k+1}=Y^ka^k$ 
in the space $\uqk\subset (\uqq)^{\ot 2}$. It is easy to see that
the operator $\ork(v)$ is diagonal in this base: 
$\ork(v)=diag(v_1,\,v_2,\,...,
\,v_{k+1})$. As for the operator $\ork(u)$, it is over-diagonal.

Using the relation $q^2\ork(u)\ork(v)-\ork(v)\ork(u)=-h_k\ork(u)$
we get $h_k=v_1-q^2 v_2$. Let us compute $v_1$ and $v_2$:
$\ork(v)f_1=\qq f_1$, i.e., $v_1=\qq$ and
$$ \ork(v)f_2=\qq P_k(a^{k-1}b+...+q^{k-2}aba^{k-2})-qP_k(q^{k-1}ba^{k-1})=$$
$$\qq P_kf_2-(q^k+q^{k-2})P_kba^{k-1}=\qq f_2-(q^k+q^{k-2})q^{k-1}\kappa f_2$$
where $\kappa=(1+q^2+...+q^{2(k-1)})^{-1}$. We use here the relation
$P_kba^{k-1}=q^{k-1}\kappa f_2$ which can be deduced from the following
equalities 
$$P_k f_2=f_2,\, P_k(a^i(qab-ba)a^{k-1-i})=0,\,\, i=0,...,k-2.$$

This gives $v_2=\qq-(q^{2k-1}+q^{2k-3})\kappa$. Finally we have
$h_k=\qq-q+q^2(q^{2k-1}+q^{2k-3})\kappa=(\qq+q^{2k+1})\kappa$. Hence, we have
proved the following
\begin{proposition}
If $\rho$ is the spin 1/2 representation of the braided 
Lie algebra $\osl$ then
the map 
$$\rho_k=(\qq+q^3)(1+q^2+...+q^{2(k-1)})(\qq+q^{2k+1})^{-1}\rho$$
is also a representation.
\end{proposition}

Thus, the quantity $(\qq+q^3)(1+q^2+...+q^{2(k-1)})(\qq+q^{2k+1})^{-1}$ is a
proper "q-analogue" of the factor $k$ in the truncated version (3) of the
tensor power.

\section{Reflection equations and related objects}

Let $\gggg$ be a classical simple Lie algebra  and ${\cal R}$ be the
quantum universal  R-matrix corresponding to $\gggg$. 
The reflection equations
(RE) are the relations $$SL_1SL_1=L_1SL_1S$$
where $L_1=L\ot id,\,L=(l_i^j),\, 1\leq i,j\leq n,\, 
S=\sigma \rho^{\ot 2}({\cal
R}),\, \rho:\gggg\to End(V),\,\,(dim\,V=n)$ is the vector fundamental
representation and  $\sigma$ is the flip. (Let us remark that there exist 
other forms of the RE, see, e.g., \cite{KSS}.)

Such equations have been introduced by I.~Cherednik for quantum R-matrices
depending on a spectral parameter. They also occur  in the construction of 
braided groups by S.~Majid.  (It is worth noting that the Majid's braided 
groups possess
a braided Hopf structure which may be introduced by
 the usual matrix coproduct and
a proper antipod, but we do not need this structure.) 

Let  $\aq$ denote the graded quadratic algebra generated by the RE
(the RE algebra). More precisely, we put
$\aq=T(W)/\{I_-^q\}$ where the space $W=Span(l_i^j)$ is
 generated by the matrix
elements  $l_i^j$ and $I_-^q$ is the
subspace in $\ww$ generated by the matrix elements 
of $SL_1SL_1-L_1SL_1S$. Thus,
the space $I_-^q$ is a q-analogue of the space of 
skew-symmetric tensors in
$\ww$  and the algebra $\aq$ is a q-analogue of the 
symmetric algebra $Sym(W)$ of $W$.

It is well-known  that this algebra becomes a $Fun_q(G)$-comodule
if equipped with the coproduct
$$\Delta:\aq\to \aq\ot Fun_q(G),\; l_i^j\to t_i^m\gamma(t_n^j)\ot l_m^n$$
where $\{t_i^j\}$ is the usual basis in the space $Fun_q(G)$ and  $\gamma$
is the antipode in $Fun_q(G)$ (cf. \cite{KSS}, \cite{I1} for details). More
precisely, if we equip the space $W$ and therefore the algebra $T(W)$ with a 
$Fun_q(G)$-comodule structure by means of the above coproduct then the space
$I_-^q$ becomes a sub-comodule of the comodule $W^{\ot 2}$.

\begin{remark} One often considers another subspace 
$$I_+^q\subset \ww,\,I_+^q=Span(SL_1SL_1+L_1SL_1S^{-1})$$
which is a q-deformation of the space of symmetric tensors. The corresponding
graded quadratic algebra $T(W)/\{I_+^q\}$ (together with the algebra
$\aq$) participates  in the construction of a quantum differential 
calculus in the algebra
$Fun_q(Mat(n))$ (cf. \cite{IP1}, \cite{IP2}). Namely, 
the generators  of the 
algebra $T(W)/\{I_+^q\}$ are considered as one-sided invariant 
differentials on
the quantum group. Meanwhile, the algebra $\aq$ is treated as one of 
invariant vector fields over quantum linear group. \end{remark}

A natural question arises: whether the quadratic graded algebra $\aq$ 
is a flat
deformation of the ordinary symmetric algebra $Sym(W)=Fun(Mat(n))$? 
The answer is negative for simple Lie algebras of 
$B_n,\, C_n,\,D_n$ series.
As for the series $A_n$ it seems very plausible that the answer is positive. 
As far as we know there are no proofs of the flatness up to now, but
a necessary condition, namely, the existence of the 
corresponding Poisson bracket, is
fulfilled.

The mentioned bracket is well defined in the space $Fun(Mat(n))$ and
it is quadratic with respect to the basis formed by the 
matrix elements $l_i^j\in
Fun(Mat(n))$ (for details we refer the reader to the paper \cite{I2}).
Let us only remark that the linearization of this bracket gives
rise to another
Poisson  bracket which is linear in generators $l_i^j$. Moreover, these two
brackets are compatible, i.e., any their linear combination 
is a Poisson bracket
as well. The family of these linear combinations is 
called a {\em Poisson pencil}.

Assuming the deformation $Sym(W)\to \aq$ to be flat (in what follows we
restrict ourselves to the $sl(n)$-case) we can consider 
the algebra $\aq$ as a
quantization of the above mentioned quadratic Poisson bracket.

It is not difficult to construct (under the above assumption of flatness)  a
quantum counterpart for the whole of the Poisson pencil. It can be done by 
a simple change of basis in the algebra $\aq$ (in \cite{GR} 
this method has been
used to quantize another Poisson pencil). 

To construct this quantum counterpart introduce the shift operator
$Sh(l_i^j)=l_i^j+h\delta_i^j$. Extending this operator to the 
algebra $\aq$ by
multiplicativity, we get a two-parameter algebra $\ahq=T(W)/\{J\}$ where the
subspace  $J\subset W\oplus \ww$ is generated by the
matrix elements of $SL_1SL_1-L_1SL_1S-h(L_1S^2-S^2L_1)$.

The latter algebra can be treated as the enveloping algebra of a
Lie algebra-like object. To do this we introduce the following 
"generalized Lie
bracket" 
$$[\,\,,\,\,]:\,I_-^q\to W,\;
[\,\,,\,\,](SL_1SL_1-L_1SL_1S)=h(L_1S^2-S^2L_1)$$
(the latter operator is considered on each matrix element separately).

The following proposition clarifies what property of the ordinary 
Lie bracket is 
inherited by the above one.

\begin{proposition} The data $(W,\,I_-^q,\,[\,\,,\,\,])$  satisfy the 
following relations \begin{equation}
([\,\, ,\,\,]^{12}-[\,\, ,\,\,]^{23})(I_-^q\ot W\cap W\ot 
I_-^q)\subset I_-^q;
\end{equation}
 \begin{equation}
[\,\, ,\,\,]([\,\, ,\,\,]^{12}-[\,\, ,\,\,]^{23})
(I_-^q\ot W\cap W\ot I_-^q)=0.
\end{equation}
\end{proposition}
{\bf Proof.}  Consider an element $Q$ of the space 
$I_-^q\otimes W\cap W\otimes I_-^q$. By definition it can be represented as 
$Q=\sum f_i k_i=\sum m_j g_j$ where $f_i, g_i$ are certain elements from 
$I_-^q$ and $k_i, m_j$ belong to $W$.
Apply the operator $[\,\,,\,\,]^{12}-[\,\, ,\,\,]^{23}$ to the element
$Q$. We have 
$$([\,\,,\,\,]^{12}-[\,\, ,\,\,]^{23})Q=\sum ([f_i]
k_i-m_i[g_i])=  \sum (f_i c_i-d_i g_i)
$$
where the elements  $c_i,\,d_i \in k$ are defined as follows: $c_i= (Sh-id)
l_i,\, d_i=(Sh-id) m_i$ (to get the latter relation it suffices to apply the
operator $Sh$ to  $0=Q-Q=\sum f_i k_i-\sum m_j g_j$ and to gather 
the elements
containing $h$). This proves the relation (6). To get (7) it suffices to 
use the
relation $\sum[f_i] c_i=\sum d_i[g_i]$ (we can prove it applying the operator
$Sh$ to $Q-Q$ and gathering the terms with $h^2$).

\begin{remark} Let us remark that the  formulae (6) and (7) represent 
the most
general form of Jacobi identity connected to the deformation theory. It was
suggested by  A.~Polishchuk and L.~Positcelski
 who have proved a version of the PBW theorem (under the assumption of
"Koszulity" of the corresponding graded quadratic algebra) 
for  Lie algebra-like objects of such type (cf. \cite{BG} where a 
slight generalization of the
mentioned result is given).

Using this result it is possible to prove the following statement. 
Let a family
of graded quadratic algebras $\aq=T(V)/\{I^q\}$ be a 
flat deformation of the
symmetric algebra $Sym(V)$ of a linear space $V$ and let a bracket
$[\,\,,\,\,]:I^q\to V$ satisfy the above version of Jacobi identity 
with
$I_-^q=I^q$. Then the corresponding enveloping algebra  (with a
factor $h$ at the bracket) $\ahq=T(V)/\{I^q-h\br I^q\}$ is also a 
flat two parameter
deformation. 

Let us also remark that the above bracket $[\,\,,\,\,]$ is defined
only on the
space $I_-^q$ but we can extend it   to the whole $W$ in a natural way
by setting $[\,\,,\,\,]\,I_+^q=0$. The extended data
$(V,\,I_-^q,\,I_+^q,\,[\,\,,\,\,])$ of such type (with a complementary
condition that the quadratic algebra $T(V)/\{I^q\}$ is Koszul)
 have been called in \cite{DG2}  generalized Lie algebras.

Nevertheless, in the case under consideration we do not need the above
mentioned generalization of the PBW theorem since the
deformation $\aq \to \ahq$ is flat by construction and 
assuming the deformation
$Sym(W) \to \aq$ to be flat we  automatically have a flat 
two parameter deformation $Sym(W) \to \ahq$. 

Let us also remark that the Jacobi identity (6), (7) differs  from
that given by the formula (5). If the former one is connected to the deformation
theory of quadratic algebras the latter one is rather motivated by 
representation theory of the braided Lie algebras.
\end{remark}

Thus, we have defined generalized Lie algebras related to the RE 
(assuming the flatness conjecture to be true, it is easy to show that the
graded quadratic algebra $\aq$ is Koszul for generic $q$).  There exists a 
number
of papers where  these  Lie algebra-like objects are considered as proper
q-deformations of  ordinary Lie algebras, taking the following
fact into consideration. The above mentioned linear Poisson bracket 
coincides  with the $gl(n)$ Lie bracket on the generators
$l_i^j$ 
(cf. \cite{KSS}, \cite{I2}). This is also a
reason why the elements of the algebra $\aq$ are often treated as one-sided
invariant vector fields.

However, in contrast to the braided Lie algebras considered in the 
previous Sections 
these
objects are constructed in a  trivial (from the deformation theory point of
view) way and for this reason they do not give rise to any interesting 
generalization of the Lie algebra structure.

Let us remark that the space $W\,\,(dim\, W=n^2)$ can be decomposed into a
direct sum of two irreducible $Fun_q(SL(n))$-comodules of dimensions 1 and 
$n^2-1$
respectively. This decomposition coincides in fact with that described in 
Section 2 (note that there exists a natural way to equip any
$Fun_q(SL(n))$-comodule with a $U_q(sl(n))$-module structure).

We can treat the latter component as a q-deformation of $sl(n)$. However, the
algebras $\ahq$ (even if $h=1$) cannot be restricted to this component (to see
this, it suffices to show that the above mentioned quadratic Poisson bracket
cannot be restricted to $sl(n)$). So, if we want to get a Lie algebra-like
object which is well-defined on the latter component, we come to the above
construction of a braided Lie bracket.

{\bf Acknowledgment} The author thank G.~Arutyunov, J.~Donin, A.~Isaev, 
P.~Pyatov, S.~Majid and S.~Shnider for helpful discussions. The work was 
partially supported by grant ISF MBI000.

 \end{document}